\documentclass{proc-l}
\usepackage{amssymb,amscd,verbatim,amsmath,color,epsfig}
\newtheorem{theorem}{Theorem}
\newtheorem{lemma}[theorem]{Lemma}
\newtheorem{cor}[theorem]{Corollary}

\newtheorem{proposition}[theorem]{Proposition}
\theoremstyle{remark}

\theoremstyle{definition}

\theoremstyle{definition}

\theoremstyle{definition}

\let\a\alpha  \let\b\beta    
  \let\g\gamma
  \let\l\lambda  \let\k\chi

\def\G{\Gamma}

\def\k{\kappa}
\def\t{\theta}

\def\Z{{\mathbb Z}}

\def\H{{\mathbb H}}

\def\SL{\textup{SL}}

\newcommand{\hp}{\mathbb{H}}

\numberwithin{theorem}{section} \numberwithin{equation}{section}

\begin{document}

\title{zeros of some level 2 Eisenstein series}

\author[Garthwaite, Long, Swisher,
Treneer]{Sharon Garthwaite, Ling Long, Holly Swisher,
Stephanie Treneer}

\address{ Bucknell University, Lewisburg, PA 17837\endgraf
sharon.garthwaite@bucknell.edu\endgraf
\null
 Department of
Mathematics, Iowa State University, Ames, IA 50011 \endgraf
linglong@iastate.edu\endgraf
\null
Oregon State University, Corvallis, OR 97301\endgraf
swisherh@math.oregonstate.edu\endgraf
\null
Western Washington University, Bellingham, WA 98225\endgraf
stephanie.treneer@wwu.edu\endgraf
}

\date{\today}
\thanks{The second author was supported in part
by the NSA grant \#H98230-08-1-0076. }

\begin{abstract}
The zeros of classical Eisenstein series satisfy many intriguing
properties.  Work of F. Rankin and Swinnerton-Dyer pinpoints their
location to a certain arc of the fundamental domain, and recent work
by Nozaki explores their interlacing property.  In this paper we
extend these distribution properties to a particular family of
Eisenstein series on $\Gamma(2)$  because of its elegant connection
to  a classical Jacobi elliptic function $cn(u)$ which satisfies a
differential equation (see formula \eqref{eq:3}). As part of this
study we recursively define a sequence of polynomials from the
differential equation mentioned above that allow us to calculate
zeros of these Eisenstein series.   We end with a result linking the
zeros of these Eisenstein series to an
 $L$-series. \end{abstract}

\maketitle
\noindent
\emph{Accepted for publication by the Proceedings of the American Mathematical Society}

\section{Introduction}

Modular forms are functions defined on the complex upper half-plane,
$\H$, that transform in a nice way with respect to subgroups of
$\mathrm{SL}_2(\Z)$. These functions have proved to be intricately
related to deep results in many areas of number theory,
combinatorics, algebraic geometry, and other areas of mathematics
(cf. \cite{ono-book-modularity}). For example,
modular forms are deeply involved in the proof of Fermat's last
theorem, results in partition theory, and denominator formulas for
the monster group. Eisenstein series are primary components of
modular forms.  That is, the weighted algebra of all integral weight
holomorphic modular forms for $\mathrm{SL_2}(\Z)$ is generated by
$E_4(z)$ and $E_6(z)$, the unique normalized Eisenstein series of
weight 4 and 6 for $\mathrm{SL}_2(\Z)$, respectively.  A nice
generating function of the classical Eisenstein series $E_{2k}(z)$
is the Weierstrass $\wp$ function, which is an elliptic function
satisfying
\begin{equation*}\label{eq:1}
  (\wp'(u))^2=4\wp(u)^3-a_4\wp(u)-a_6,
\end{equation*}
where $a_4,a_6$ are scalar multiples of $E_4$ and $E_6$,
respectively.

Location of the zeros of Eisenstein series is of
fundamental importance. Provided with such information,  Li, Long,
and Yang construct noncongruence  cuspforms from congruence
Eisenstein series \cite{lly05}.  Additionally, in a paper of Ono
and Papanikolas \cite{Ono-Papanikolas-04}, an interesting
formula is highlighted which relates the location of the zeros
of $E_{2k}(z)$ to special values of the Riemann Zeta function
$\zeta(s)$:
\begin{equation}\label{eq:2}
\frac{2}{\zeta(1-2k)}=60\cdot 2k-\sum_{\tau\in \H /SL_2(\Z)}
e_{\tau} \text{ord}_{\tau} (E_{2k}(\tau)) j(\tau),\end{equation}
where $e_{i}=1/2, e_{\frac{-1+\sqrt{-3}}{2}}=1/3, $ and $e_{\tau}=1$
for non-elliptic points, $\tau$.

In \cite{Rankin-SD70ZeroofEisen}, F.K.C. Rankin and H.P.F.
Swinnerton-Dyer  study  the location of the zeros of  $E_{2k}(z)$.
By using an elementary and elegant argument, they demonstrate that
all zeros of $E_{2k}(z)$ inside the fundamental domain
for $SL_2(\Z)$  are located on the arc $\{z=e^{ i \theta},
\pi/2\le \theta\le 2\pi/3\}$. Equivalently, the $j$-values of these
zeros are all real and are within the interval $[0,1728]$, where
$j(z)$ is the well-known classical modular function which
parametrizes isomorphism classes of elliptic curves. In
\cite{koh04-Eisen}, Kohnen gives an explicit formula for the zeros
of $E_{2k}(z)$.

Recently, there are several results generalizing
\cite{Rankin-SD70ZeroofEisen} to Eisenstein series of other groups
\cite[etc.]{Hah07,MNS-Eisen23-07,Shigazumi07}, other modular forms
\cite{Getz04, Gun06},  and certain weakly
holomorphic modular forms  \cite{D-J08}. Moreover, Nozaki proves
that the zeros of the classical Eisenstein series $E_{2k}(z)$ interlace
with the zeros of $E_{2k+12}(z)$ \cite{Nozaki08}.  Despite these
results, many mysteries remain about Eisenstein series, and their
zeros. For example, Gekeler studies polynomials $\varphi_{2k}(x)$
which encode the $j$-values of the non-elliptic zeros of $E_{2k}(z)$
\cite{Gekeler01Eisenstein}. He conjectures that the
$\varphi_{2k}(x)$ are irreducible with Galois group $S_d$ where $d$
is the degree of $\varphi_{2k}(x)$.\medskip

In Section \ref{notation} of this paper, we define a series of odd
weight Eisenstein series denoted by $G_{2k+1}(z)$ for $\Gamma(2)$,
the principal level 2 congruence subgroup.
 In
\cite{ly05}, Long and Yang use $G_{2k+1}(z)$ to give a second proof
to some beautiful formulae of Milne \cite{milne-square} on
representing natural numbers in terms of sums of squares or
triangular numbers. Our attention was drawn to this
family of Eisenstein series $G_{2k+1}(z)$  largely because it has an
elegant generating function $cn(u)$, a classical Jacobi
elliptic function \cite{Hancock} satisfying
\begin{equation}
\label{eq:3}\left(\frac{\textrm{d} \, cn(u)}{\textrm{d}\, u}
\right)^2 =(1-cn^2(u))\left(1-\lambda+\lambda
cn^2(u)\right).\end{equation} Here $\lambda$ is the classical
lambda function which parameterizes all isomorphism classes of
elliptic curves with full 2-torsion structure.  The equation
\eqref{eq:3} gives rise to recursions satisfied by $G_{2k+1}(z)$.
Consequently, one can compute the $\lambda$-values of the zeros of
$G_{2k+1}(z)$  with great convenience and ease.  Data from the first
few dozen  $k$ values indicates that all $\lambda$-values of the
zeros of $G_{2k+1}(z)$ are real and are within the interval
$(-\infty,0]$. Moreover, this data indicates that the zeros of
$G_{2k-1}(z)$ interlace with the zeros of $G_{2k+1}(z)$.

In Section \ref{sec:3} we investigate the location of the zeros of
$G_{2k+1}(z)$ within the following fixed fundamental domain of
$\G(2)$ (cf. Figure 1 in Section \ref{notation}):
\[
D=\left\{z\in \H: \ -1\le \text{Re}(z)\le 1, \ |z- 1/2|\ge1/2, \  \
|z+1/2|\ge 1/2\right\}.
\]

\noindent A direct application of the argument of Rankin and
Swinnerton-Dyer shows that at least one third of the zeros
of $G_{2k+1}(z)$ lie in \[\{z \in \H: \ |z+1/2|=1/2\},\] or
equivalently, have $\l$-values in $(-\infty,0]$. By refining this
argument, we will show  the following improvement.
\begin{theorem}\label{thm:main}
At least  90\% of  the zeros of
$G_{2k+1}(z)$ have real $\l$-values in the range $(-\infty,
0]$.
\end{theorem}

 In addition, by restricting our domain slightly we
obtain the following result about the separation property of the
zeros of $G_{2k-1}(z)$ and $G_{2k+1}(z)$ following an approach of
Nozaki \cite{Nozaki08}. We state this result in terms of the related
function $F_{2k+1}(z_\t)$ defined in $\eqref{EQ}$ where
$z_\t:=\frac{1}{2}e^{i\t} -\frac{1}{2}$.  Before stating the result,
we define for an integer $k>15$ the intervals
\[
I_{j,2k-1}=\left (\a_{j,k}-\frac{2\pi}{(2k+1)(2k-1)},
\a_{j,k}+\frac{2\pi}{(2k+1)(2k-1)}\right )
\]
for each $j=0,\ldots,k-1$, where $\a_{j,k}:=2\pi j/(2k-1)$.  Furthermore, let
\[
I_{2k-1}:= \bigcup_j I_{j,2k-1}.
\]
\begin{theorem}
\label{interlacing} Let $k>15$ be an integer, and $I_{j,2k-1}$,
$I_{2k-1}$ as defined above. Then the zeros of $F_{2k-1}(z_\t)$ and
$F_{2k+1}(z_\t)$ in $[\pi/10,9\pi/10]$ are restricted to $I_{2k-1}$
and $I_{2k+1}$ respectively (in fact each $I_{j,2k-1}$ and
$I_{j,2k+1}$ in $[\pi/10,9\pi/10]$ contains an odd number of zeros).
Moreover, the intervals are pairwise disjoint, and the $I_{j,2k-1}$
are separated by the $I_{j,2k+1}$.
\end{theorem}

In Section \ref{L}, we obtain a formula (Corollary
 \ref{cor:3.2}) similar to equation \eqref{eq:2} which relates the sum of the
$\lambda$-values of the zeros of $G_{2k+1}(z)$ to the special values
of an explicit L-series.

We note that the properties of our Eisenstein series $G_{2k+1}(z)$
are similar to other families of level 2 Eisenstein series studied
by Tim Huber to which we expect our results can be extended.
Furthermore, the behavior of our Eisenstein series is
similar to all other known cases for congruence subgroups of genus
zero with relatively simple fundamental domains. It would be very
interesting to know what to expect for Eisenstein series or cuspforms
that are invariant under congruence subgroups with more complex
fundamental domains or higher genus.

The authors would like to acknowledge that this project and
collaboration was initiated at the Women in Numbers Conference at
BIRS.  The authors thank organizers Kristen Lauter, Rachel Pries and
Renate Scheidler, along with BIRS, the Fields Institute, Microsoft
Research, the NSA, PIMS, and the University of Calgary for their
generous support of this conference.

\section{ $\G(2)$ and its modular forms}

\label{notation}

Let $\G(2)$ be the principal level 2 congruence subgroup consisting
of matrices which become the identity under the natural modulo 2
homomorphism. This group is of genus 0 and has a fundamental domain
shown in Figure 1.
\begin{figure}
\includegraphics*[scale=.6]{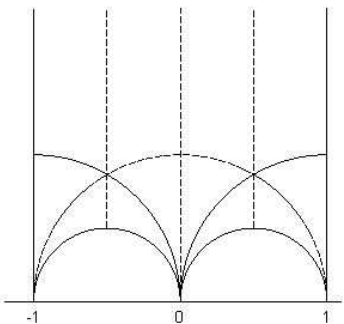}  \caption{A fundamental
domain for $\Gamma(2)$} \end{figure} Let $q=e^{2\pi iz}$; the Jacobi
theta functions are defined by
\[
 \Theta_2(z):=q^{1/8}\sum_{n\in\Z}q^{n(n+1)/2},\;\Theta_3(z):=\sum_{n\in\Z}q^{n^2/2},\!\mbox{ and }\Theta_4(z):=\sum_{n\in\Z}(-1)^nq^{n^2/2}.
\]
The square of each of the above series is a weight 1 modular form
for $\G(2)$. The classical lambda function is defined as
$\displaystyle \lambda(z):=\frac{\Theta_2^4(z)}{\Theta_3^4(z)}$ and
it generates the field of meromorphic modular functions for $\G(2)$.

\subsection{Some odd weight level 2 Eisenstein series}

We now focus on a particular family of Eisenstein series of
odd weight on $\G(2)$ investigated in \cite{lly05}.
To begin, we define the character
\[ \chi(\g) := \chi(d) =
\left( \frac{-1}{d}\right)\] where $\gamma=\left(\begin{smallmatrix}
a & b\\ c & d \end{smallmatrix} \right)$ and $ \left (
\frac{-1}{\bullet} \right )$ denotes a Legendre symbol. We now
define
\begin{equation*}\label{eq:E2k+1}
E_{2k+1,\chi}(z) := \frac{1}{2}\sum_{\substack{(c,d)\equiv (0,1) \;
(2)\\ (c,d)=1}}\chi(d)\cdot {(cz+d)^{-(2k+1)}},
\end{equation*}
and recall some results about $E_{2k+1,\chi}(z)$ (see \cite[Lemmas
5.2 - 5.5.]{lly05}).

\begin{lemma}
\label{FourierEx} Define the $n$th Euler number $e_n$ by $
\displaystyle \sec(t) = \sum_{n=0}^\infty \frac{e_n}{n!}t^n. $ The
function $E_{2k+1,\chi}(z)$ is a modular form of weight $2k+1$ on
$\G(2)$ with Fourier expansion at $\infty$ given by
\[ E_{2k+1,\chi}(z) = 1
+ \frac{4(-1)^k}{e_{2k}}\sum_{r=1}^\infty
\chi(r)\frac{r^{2k}q^{r/2}}{1-q^{r/2}}. \] \end{lemma}

For any integer $k$ let
\begin{equation}\label{Sk}
S_k(\a,\b)=S_k(\a,\b)(z):=\sum_{\substack{(c,d) \equiv (\a,\b) \;
(4)\\ (c,d)=1}} \frac{1}{(cz+d)^k}.
\end{equation}
By the explicit description of $\chi$,
\[
E_{2k+1,\chi}(z) = \frac{1}{2}\left( S_{2k+1}(0,1) + S_{2k+1}(2,1) -
S_{2k+1}(0,3) - S_{2k+1}(2,3)\right).
\]

Recall that if $f$ is a meromorphic function on $\H$,
$\gamma=\left(\begin{smallmatrix} a & b\\ c & d
\end{smallmatrix} \right)\in\mathrm{SL}_2(\Z)$, and $k$ is an integer
then the slash operator in weight $k$ is
\[f(z)|_k \g :=(cz+d)^{-k}f(\g z).\]
We use $f|_\g$ to denote $f(z)|_k~\g$. The following lemma for odd $k$
follows analogously from Lemma 3.3 in \cite{lly05}.
\begin{lemma}\label{S-lemma}
If $\gamma\in \SL_2(\Z)$, and $S_k$ is defined as above for odd $k$,
then
\begin{equation*}\label{eq:S_k}
S_k(\a,\b)\mid_\g = S_k((\a,\b)\cdot\g) \pmod{4})
\end{equation*}
where $(\a,\b)\cdot\g \pmod{4}$ denotes matrix multiplication modulo
4. \end{lemma}

\subsection{A generating function}\label{subsec:gen.fcnfor E2k+1}
Recall that $cn(u)$ is a classical Jacobi
elliptic function which satisfies (cf. \cite[pp. 256]{Hancock})
\begin{equation*}
\frac{\k K}{2\pi}cn\left(\frac{2Ku}{\pi}\right)=\frac{\sqrt{q}\cos
u}{1+q}+\frac{\sqrt{q^3}\cos 3u}{1+q^3}+\frac{\sqrt{q^5}\cos
5u}{1+q^5}+\cdots,
\end{equation*}
where $\sqrt{\k}=\displaystyle\frac{\Theta_2(2z)}{\Theta_3(2z)}$ (cf
\cite[pp. 241]{Hancock}) and hence $\l(2z)=\k^2$ and
$K=\frac{\pi}{2}\Theta_3^2(2z)$.  For the remainder of this section
we simply write $\l$ for $\l(2z)$.  To see how the family of
Eisenstein series $E_{2k+1,\chi}$ relates to $cn(u)$,  we consider
the expansion of $E_{2k+1,\chi}$ at the cusp $1$, given by $E_{2k+1,\chi}|_{\g_0}$, where
$\g_0=\left(\begin{smallmatrix} 1&-1\\1&0
\end{smallmatrix}\right)$. The previous lemma allows us to easily
calculate $E_{2k+1,\chi}\mid_{\g_0}$ in terms of the $S_k$, which
brings us to the following conclusion.
\begin{align*}
E_{2k+1,\chi}\mid_{\g_0}(z)  = \frac{1}{2}\left( S_{2k+1}(1,0) +
S_{2k+1}(3,2) - S_{2k+1}(3,0) - S_{2k+1}(1,2)\right).
\end{align*}

It is this function $G_{2k+1}(z):=E_{2k+1,\chi}\mid_{\g_0}(z)$ whose
zeros we study. It has been shown (\cite{ly05} Lemma 5.5) that
\begin{equation*}
G_{2k+1}(z)=\frac{4(-i)^{2k+1}}{e_{2k}}\sum_{r=1}^{\infty}\frac{(2r-1)^{2k}q^{(2r-1)/4}}{1+q^{(2r-1)/2}}.
\end{equation*}
In particular, $G_1(z)=-i\Theta_2^2(z)$ (cf. \cite{ly05} Lemma
5.4).

 Using the series
expansion $\displaystyle
  \cos(x)=\sum_{n=0}^{\infty}(-1)^n\frac{x^{2n}}{(2n)!}, $
we have
\begin{eqnarray*}
\frac{\k K}{2\pi}cn\left(\frac{2Ku}{\pi}\right)&=&\sum_{k=0}^{\infty}\frac{(-1)^k}{(2k)!}\left
(\sum_{r=1}^{\infty}\frac{(2r-1)^{2k}q^{(2r-1)/2}}{1+q^{2r-1}}\right)u^{2k}\\ &=&\sum_{k=0}^{\infty}
\frac{ie_{2k}}{4(2k)!}G_{2k+1}(2z) u^{2k}.
\end{eqnarray*}
As an immediate consequence, we have the following proposition.
\begin{proposition}
Let $\k = \frac{\Theta_2^2(2z)}{\Theta_3^2(2z)}$ and
$K=\frac{\pi}{2}\Theta_3^2(2z)$.  We have that
\begin{equation}\label{eq:3.4}
 cn(u)=\sum_{k=0}^{\infty}
  \frac{ie_{2k}\pi^{2k+1}}{2^{2k+1}(2k)!}\frac{G_{2k+1}(2z)}{K^{2k+1}\k}u^{2k}.
\end{equation}
\end{proposition}

From  $cn (u)$ satisfying \eqref{eq:3} (cf. \cite[pp.
247]{Hancock}), we have that
\begin{equation}\label{eq:3.5}
  cn(u)=1-\frac{u^2}{2!}+(1+4\l)\frac{u^4}{4!}
  -(1+44\l+16\l^2)\frac{u^6}{6!}+\cdots.
\end{equation}
Note that the choice $K=\frac{\pi}{2}\Theta_3^2(2z)$ aligns the
constant terms of equations \eqref{eq:3.4} and \eqref{eq:3.5}.
We have that  \[ \frac{(-1)^ki\pi^{2k+1}e_{2k}}{2^{2k+1}} \cdot
\frac{G_{2k+1}(2z)}{K^{2k+1}\k}
=\frac{(-1)^k{i}e_{2k}G_{2k+1}(2z)}{\Theta_3^{4k}(2z)\Theta_2
^2(2z)}=
\frac{(-1)^ke_{2k}G_{2k+1}(2z)}{\Theta_3^{4k}(2z)G_1(2z)}\] is a
modular function for $\G(2)$ which is holomorphic everywhere
except possibly a pole of finite order at the cusp infinity
coming from $\Theta_3^{4k}(2z)$. Consequently, for $k\geq 1$ it is a degree
$(k-1)$ polynomial in terms of $\l$, denoted $p_{2k+1}(\l)$, having the same zeros as $G_{2k+1}$ in the fundamental domain.  We list the first six below.  \\

\begin{tabular}{ll}
$p_1(\l)=1$, & $p_7(\l)=1+44\l+16\l^2$,\\
$p_3(\l)=1$, & $p_9(\l)=1+408\l+912\l^2+64\l^3$,\\
$p_5(\l)=1+4\l$, &
$p_{11}(\l)=1+3688\l+30764\l^2+15808\l^3+256\l^4$. \\
& \\
\end{tabular}

There are $k-1$  nontrivial zeros counting multiplicity, in addition
to the trivial zeros of $G_{2k+1}$ at the cusps. Calculations for
$2k+1\leq 51$ reveal some striking numerical trends. The $\l$-values
of the zeros of $G_{2k+1}$ all lie in $(-\infty, 0)$.  We list the
numerical zeros of the first few $p_{2k+1}(\l)$ below.  \\\\
\noindent
\begin{tabular}{|l|c||l|c|} \hline   k=2& \{$-0.25$\} &
k=4&\{$-0.0025, -0.4598, -13.788$\}\\ \hline k=3&\{$-.0229,
-2.7271$\} & k=5&\{$-.00027, -.1280, -1.8792, -59.7425$\}\\
\hline \end{tabular}\\\\ Moreover, the $\l$-zeros of $p_{2k-1}(\l)$
interlace with the $\l$-zeros of $p_{2k+1}(\l)$. In addition, the
polynomials $p_{2k+1}(\l)$ are irreducible with Galois group
$S_{k-1}$ for $k\le 9$.

\section{Locating the zeros of $G_{2k+1}(z)$}\label{sec:3}

Based on the numerical evidence above, we now turn our focus to the
$\l$-values of the zeros of $G_{2k+1}(z)$.  Proving that the
$\l$-values lie in $(-\infty,0]$ is equivalent to proving that the
zeros of $G_{2k+1}(z)$ lie on the line Re$(z)=1$, a boundary arc of the
fundamental domain for $\G(2)$. This, in turn, is equivalent to proving that the zeros of $G_{2k+1}\mid_{\g_1}(z)$, where
$\g_1=\bigl(
\begin{smallmatrix} 0 & -1\\ 1 & 0
\end{smallmatrix} \bigr)$, lie on the arc $|z+\frac{1}{2}|=\frac{1}{2}$, as we mentioned in Section 1.   Our main strategy is
to generalize the approach of Rankin and Swinnerton-Dyer.  In
particular, we will need to refine estimates of the terms arising in
this approximation.

We first find, using Lemma \ref{S-lemma}, that
\begin{equation*}
G_{2k+1}\mid_{\g_1}(z) = \frac{1}{2}\left(S_{2k+1}(0,3) +
S_{2k+1}(2,1) -S_{2k+1}(0,1) - S_{2k+1}(2,3)\right).
\end{equation*}

We now wish to count the number of zeros of $G_{2k+1}\mid_{\g_1}(z)$ that must lie on the arc $|z+\frac{1}{2}|=\frac{1}{2}$.

\begin{lemma}\label{newS} Let $\alpha$ be even and $z_\t=\frac{1}{2}e^{i\t} -\frac{1}{2}$.  Then
\[
S_k(\alpha,\beta)(z_\t)=S_k\left(\frac{\alpha}{2},\beta-\frac{\alpha}{2}\right)(e^{i\t})+S_k\left(\frac{\a}{2}+2,\b-2-\frac{\a}{2}\right)(e^{i\t}).
\]
\end{lemma}

\begin{proof} The lemma follows by separating the terms in $S_k(\a,\b)(z_\t)$ for pairs $(c,d)$ according to the two possible values of $\frac{c}{2}$ modulo $4$, and by the fact that $(c,d)=1$ implies $(c,d)=(\frac{c}{2},d-\frac{c}{2})$.\end{proof}

By Lemma \ref{newS} and the observation that $S_k(-a,-b)=-S_k(a,b)$,
we have
\[
G_{2k+1}\mid_{\g_1}(z_\t)=\left(S_{2k+1}(0,3)+S_{2k+1}(1,0)+S_{2k+1}(2,1)+S_{2k+1}(3,2)\right)(e^{i\theta}).
\]
We next factor out $(e^{\frac{i\theta}{2}})^{-2k-1}$ following
\cite{Rankin-SD70ZeroofEisen} and define
\begin{multline}
\label{EQ}
F_{2k+1}(z_\t):=(e^{\frac{-i\t}{2}})^{-2k-1}G_{2k+1}\mid_{\g_1}(z_\t)\\
= \left(\widetilde{S}_{2k+1}(0,3)+\widetilde{S}_{2k+1}(1,0)+\widetilde{S}_{2k+1}(2,1)+
\widetilde{S}_{2k+1}(3,2)\right)(\t),
\end{multline}
where
\[\widetilde{S}_k(\alpha,\beta)(\theta):=\sum_{\substack{(c,d) \equiv (\a,\b) \;(4)\\ (c,d)=1}}(ce^{i\t/2}+de^{-i\t/2})^{-k}.\]

From equation \eqref{EQ} we see that $F_{2k+1}(z_\theta)$ is purely
imaginary by noting that
$F_{2k+1}(\overline{z}_\theta)=-F_{2k+1}(z_\theta)$, as conjugation
interchanges each pair of sums. Moreover, in each sum
$c$ and $d$ have opposite parity.

A direct application of the analysis in
\cite{Rankin-SD70ZeroofEisen} shows there are exactly the right number
of zeros when $\t$ ranges through $[\frac{\pi}{3}, \frac{2\pi}{3}]$. In
other words, one third of the zeros are expected to be on the
above arc with $\t$ in $[\frac{\pi}{3}, \frac{2\pi}{3}]$. The real
difficulty comes from the analysis when $\t$ is close to 0 or $\pi$.

\subsection{Extraction of the main term}

To begin our analysis of $G_{2k+1}\mid_{\g_1}(z_\t)$ on the arc $| z
+ \frac{1}{2}| =\frac{1}{2}$,  we first extract the two terms for which $c^2+d^2 = 1$ (one occurring in each of $\widetilde{S}_{2k+1}(0,3)$ and
$\widetilde{S}_{2k+1}(1,0)$)  to create our
main term, \[
(e^{\frac{i\t}{2}})^{-2k-1}-(e^{\frac{-i\t}{2}})^{-2k-1} =
-2i\sin\left(\frac{\t(2k+1)}{2}\right). \] Thus, we have
\begin{equation}\label{eq:F_2k+1}
F_{2k+1}(z_\t) = -2i\sin\left(\frac{\t(2k+1)}{2}\right)
 +R_{2k+1}(z_\theta),
\end{equation}
where $R_{2k+1}(z_\theta)$ or simply $R(z_\theta)$ is purely
imaginary and denotes the error term obtained by summing over all
remaining terms.

\subsection{Bounding the error term}

We now turn to bounding $R(z_\theta)$.  Here we will consider the
contribution of terms satisfying $c^2+d^2 = N$, for $N>1$.  We
take advantage of the symmetry in our remaining terms and define
the following partial sum.
For each ordered pair of nonnegative integers $(a,b)$ chosen such that $a$ is odd and $b$ is even to avoid double-counting, define
\begin{equation}
\label{P} P(a,b)(\theta)= \sum_{\substack{(|c|,|d|)=(a,b) \scriptsize \mbox{ or } \normalsize (b,a)}} (ce^{i\t/2}+de^{-i\t/2})^{-2k-1},
\end{equation}
where the sum contains only terms occurring in $F_{2k+1}(z_\t)$.  For example,
\begin{align*} P(3,0)&=
\left((-3e^{i\t/2})^{-2k-1}+ (3e^{-i\t/2})^{-2k-1}\right) =
\frac{2i}{3^{2k+1}}\sin\left(\frac{\t (2k+1)}{2}\right).
\end{align*}
Due to symmetry, each partial sum $P(a,b)$ is purely imaginary as
before.

\subsection{The case when $N\leq100$}

We now give upper bounds on the terms for which $1< N\leq 100$. We
assume that $2k+1 > 51$ as we have numerically verified the location
of zeros for low weight cases.

\noindent  When $b=0$ we have the following cases:
\[
(a,b) \in \{(3,0), (5,0), (7,0), (9,0)\}.
\]
Here,
\[
|P(a,b)(\theta)| =
\left|\frac{2}{a^{2k+1}}\sin\left(\frac{(2k+1)\t}{2}\right) \right|
\leq  \frac{2}{a^{51}}.\]
 The contribution from
these terms to the error term is smaller than
\[
E_1=\sum_{n=1}^{4} \frac{2}{(2n+1)^{51}}< 10^{-24}.
\]

\noindent  When $b\ne 0$, then for each term in $P(a,b)$ we have
\[|ce^{i\t/2}+de^{-i\t/2}|=|ae^{i\t/2}\pm be^{-i\t/2}|\mbox{ or }|be^{i\t/2}\pm ae^{-i\t/2}|.\]
In the first case,
\begin{align}\label{Pab}
|ae^{i\t/2}\pm be^{-i\t/2}|^{-2k-1}
&=\left
|(a\pm b)\cos\left(\t/2\right)+i(a\mp b)\sin\left(\t/2\right)
\right |^{-2k-1}\notag\\&=\begin{cases}
\left((a-b)^2+4ab\cos^{2}\left(\t/2\right)\right)^{-k-1/2} &
\text{if } +,\\
\left((a-b)^2+4ab\sin^{2}\left(\t/2\right)\right)^{-k-1/2} & \text{if }-.
\end{cases}
\end{align}
By symmetry, the second case yields the same result. When $a$ and
$b$ are nonzero,  $P(a,b)$ contains four terms, two of each type in
\eqref{Pab}.  Therefore,
\[
|P(a, b)|\leq \frac{2}{((a-b)^2+4ab\sin^2(\t/2))^{k}}+\frac{2}{((a-b)^2+4ab\cos^2(\t/2))^{k}}.
\]
Note that
\[
\max(\cos^2(\t/2),\sin^2(\t/2))\ge 1/2.
\]
Moreover, if we limit our choice of $\t$ to the
interval {$(0.05\pi, 0.95\pi)$} then \[
\min(\cos^2(\t/2),\sin^2(\t/2))> (.079)^2. \] Hence we conclude that
when $b$ is nonzero,
\begin{equation}\label{Pabbound}
|P(a, b)|\leq
\frac{2}{((a-b)^2+4ab({.079})^2)^{k}}+\frac{2}{((a-b)^2+2ab)^{k}}.
\end{equation}
Now consider the cases
\[
(a,b) \in \{(1,2),
(3,2), (3,4), (5,4), (5,6), (7,6)\},
\]
with $|a-b|=1$. Summing the bounds in \eqref{Pabbound} for the six
pairs $(a,b)$ listed above yields that the contribution from these
terms is smaller than
\[
E_2={0.656}.
\]

\noindent The remaining cases with $N\leq 100$ are
\[
(a,b) \!\in\!
\{(1,4), (1,6), (1,8), (3,6),(3,8), (5,2), (5,8),
(7,2),(7,4), (9,2), (9,4)\}.
\]

\noindent For these we note that $|a-b|\ge 3$. Reasoning as above
yields \[ |P(a,b)|< \frac{4}{9^{k}}.\] Hence, the total contribution
to the error term is bounded by \[E_3=11\cdot \frac{4}{9^{k}}\le
\frac{4}{9^{25}}< 10^{-10}.\]

In total, we have
\[\left|\sum_{a^2+b^2\le 100} P{(a,b)}\right|\le E_1+E_2+E_3< {0.657}.\]

\subsection{The case when $N>100$}

We now consider terms $c^2+d^2=N$ with $N\geq 101$, again
assuming that $2k+1> 51$.  The number of terms satisfying
$c^2+d^2=N$ is at most $2(2N^{1/2}+1)\leq 5N^{1/2}$. Note that
\[ |ce^{i\theta/2}+de^{-i\theta/2}|^2 = c^2+2cd\cos\theta+d^2.
\] If we restrict ourselves to $\theta$ values with
$|\cos\theta|\leq \alpha$, for some $\alpha \in (0,1)$,  then \[
c^2+2cd\cos\theta+d^2 \geq (1-\alpha)(c^2+d^2). \] Thus, we have
\[ |R(z_\theta)| < {0.657}+\sum_{N=101}^\infty
5N^{1/2}((1-\alpha)N)^{-k-1/2}. \] By bounding this latter sum with
an appropriate integral, we have \[ |R(z_\theta)|< {0.657}+
(1-\alpha)^{-k-\frac{1}{2}} \left(\frac{5}{k-1}\cdot
100^{-k+1}\right). \] We pick $\alpha = {.9877}$, which corresponds to
{$\theta \in (0.05\pi, 0.95\pi)$}. Then $|R(z_\theta)|<2$.

\subsection{Proof of Theorem \ref{thm:main}}

Following \cite{Rankin-SD70ZeroofEisen}, we consider the values
$\theta$ for which $2\sin\left(\frac{\t(2k+1)}{2}\right) = \pm 2$
and then apply the Intermediate Value Theorem to the function
$iF_{2k+1}(z_\t)$.  In doing so we conclude that this function, and
hence also $G_{2k+1}\mid_{\g_1}(z_\t)$, must have at least $1$ zero
in each interval of the form \[ \left[\frac{A\pi}{2k+1},
\frac{(A+2)\pi}{2k+1}\right]{\subset(0.05\pi,0.95\pi)}, \] where $A$
is an odd positive integer.  We also recall that $G_{2k+1}$ has
altogether $k-1$ nontrivial zeros. Thus, having previously
considered $2k+1\leq 51$ computationally, we can now verify that at
least 90\% of the zeros of $G_{2k+1}(z)$ do indeed lie on the
boundary $Re(z)=1$, and hence have real $\lambda$-values in
$(-\infty,0]$. This concludes the proof of Theorem \ref{thm:main}.

\subsection{Proof of Theorem \ref{interlacing}}

Now that we have pinpointed the location of the zeros of
$G_{2k+1}(z)$, we turn our attention to the relationship between
the zeros of $G_{2k-1}(z)$ and $G_{2k+1}(z)$.  Numerical
evidence suggests that the zeros of $F_{2k-1}(z_\t)$ and
$F_{2k+1}(z_\t)$, the shifted Eisenstein series, interlace in the same manner as demonstrated by Nozaki for the classical case with $E_{k}(z)$ and $E_{k+12}(z)$.  We now follow Nozaki's strategy by relating
the zeros of $F_{2k+1}(z_\t)$ to the zeros of the main
term of this series, a well-understood trigonometric function
with regularly spaced zeros. We then show that the additional
error terms for the Eisenstein series will not cause the zeros of $F_{2k+1}(z_\t)$
to stray far from the zeros of the main term.

We must first strengthen the bounds on the error
estimates given in the previous section. In order to do this we
will need to change our lower bound on
$\sin\left(\theta/2\right)$ and $\cos\left(\theta/2\right)$ near
the cusps $\theta =0,\pi$. Thus, we narrow our focus to theta values in the range $[\pi/10,
9\pi/10]$.

Recall that the main term of $F_{2k+1}(z_\theta)$
given in \eqref{eq:F_2k+1} is
\[f_{2k+1}(\theta) :=
-2i\sin\left(\frac{\theta(2k+1)}{2}\right).\]
  Let $\a_{j,k}$ denote the zeros of $f_{2k-1}$ on
$(0,\pi)$, that is,
\[ \a_{j,k} = \frac{2\pi j}{2k-1}, j=1,\dots, k-1.
\]

\noindent  Note that {$\a_{j,k+1}<\a_{j,k}<\a_{j+1,k+1}$ for
$j=1,...,k-1$, and}
\[
\frac{1}{2}\min_{j=1,\dots,k-1}\{\a_j-\b_j, \b_{j+1}-\a_j\} =
\frac{2\pi}{(2k-1)(2k+1)}. \]  Our aim is to show that the zeros of
$F_{2k-1}(z_\theta)$ and $F_{2k+1}(z_\theta)$ are within this
distance of $\a_{j,k}$ and $\a_{j,k+1}$, respectively.

We now revisit the error estimates for \eqref{P} specifically for the
range $\theta \in [\pi/10, 9\pi/10]$.   We start by choosing a
parameter $\gamma=0.1562$ in order that $1+8 \gamma^2>1.195$.  Note
that for $\t \in [\pi/10, 9\pi/10]$ we have both $\sin(\t/2),\cos(\t/2)\ge \gamma$.

Next note that in terms of $k$, the contribution of $|P(a,b)|$ from
the terms $\{(3,0),(5,0),(7,0),(9,0)\}$ is smaller than
$\displaystyle e_1=\sum_{n=1}^4\frac{2}{(2n+1)^{2k+1}}\le
\frac{8}{3^{2k+1}}$. We now use our new $\theta$ bounds to
strengthen the bound on the second error term arising from the set
of pairs of the form $(a,a\pm1)$ with $N\le 100$. Replacing
${0.079}$ with $\gamma$ in \eqref{Pabbound} yields
\[
 |P(1,2)| \leq \frac{2}{1.195^{k}}+\frac{2}{5^k}.
\]
The contribution from $(3,2),(3,4),(5,4),(5,6),(7,6)$ is smaller
than
\[\frac{10}{1.585^k}+\frac{10}{13^k}.
\]
Thus the total contribution from these terms is less than
\[
e_2=\frac{2}{1.195^{k}}+\frac{2}{5^k}+\frac{10}{1.58^k}+
\frac{10}{13^k}.
\]
We again note that the contribution from the remaining cases
with $N\le 100$ is $e_3\le \frac{44}{9^k}$.  Similarly, the
contribution from terms with  $N>100$, is
\begin{align*}
e_4&=\sum_{N=101}^{\infty} 5N^{1/2}((1-C)N)^{-k-1/2}<
5(1-C)^{-k-1/2}\frac{1}{k-1}100^{-k+1}\\
&<\frac{500}{(k-1)\sqrt{1-C}}(100(1-C))^{-k}<\frac{2283}{(k-1)(4.8)^k},
\end{align*}
where $C=.952$, so that $|\cos(\t)|<C$ on $(\pi/10,9\pi/10)$.

Putting the error terms altogether, we have
\[
|R_{2k+1}(z_{\t})|<g(k):= \frac{8}{3^{2k+1}}+
\frac{2}{1.195^{k}}+\frac{2}{5^k}+\frac{10}{1.585^k}+\frac{10}{13^k}+
\frac{44}{9^k}+ \frac{2283}{(k-1)(4.8)^k}.
\]

\noindent Note that the value of $g(k)$ goes to 0 rapidly when $k$ gets large.

We now return our focus to the relationship between the zeros of
$F_{2k+1}$ and the zeros of $f_{2k+1}$.  Let
\[h(x)=\frac{2\pi}{2x+1}-\frac{2\pi^3}{6(2x+1)^3} - g(x-1).\]
Note that  $\frac{2\pi^3}{6(2x+1)^3}+g(x-1)$ decays more rapidly
than $\frac{\pi}{2x+1}$; hence, we see that for $k\ge 15$ we have
$h(k)>0$. Similarly, note that as  $\sin(x)\ge x-\frac{x^3}{6}$ for
any $0\le x<1$, we can conclude that
\[
2\sin\left(\frac{\pi}{2k+1}\right)\ge
\frac{2\pi}{2k+1}-\frac{2\pi^3}{6(2k+1)^3} = h(x)+g(x-1).
\]

\noindent Consider an interval
\[I_{j,2k-1}=\left ( \a_{j,k}-\frac{2\pi}{(2k+1)(2k-1)},
\a_{j,k}+\frac{2\pi}{(2k+1)(2k-1)}\right )\cap [\pi/10, 9\pi/10]\]
about $\a_{j,k}$, a zero of $f_{2k-1}$.  We want to show that a zero
of $F_{2k-1}$ lies in this interval. By  construction,
$I_{j,2k-1}\cap I_{i,2k+1}=\emptyset$ for any integers $i,j$.

Note that
\begin{multline*}
if_{2k-1}\left(\a_{j,k}\pm
\frac{2\pi}{(2k+1)(2k-1)}\right)= \\
2\sin\left(j\pi\pm \frac{\pi}{2k+1}\right) =2(-1)^j\sin\left(\pm\frac{\pi}{2k+1}\right).
\end{multline*}
Hence {for even $j$},
\begin{align*}
iF_{2k-1}\left(\a_{j,k}+ \frac{2\pi}{(2k+1)(2k-1)}\right)
=&2\sin\left(\frac{\pi}{2k+1}\right)
+iR_{2k-1}\left(z_{\a_j+\frac{2\pi}{(2k+1)(2k-1)}}\right)\\
\geq& 2\sin\left( \frac{\pi}{2k+1}\right)-g(k-1)  \geq h(k),
\end{align*} which is positive. Similarly, $iF_{2k-1}\left(\a_j-
\frac{2\pi}{(2k+1)(2k-1)}\right){\leq -h(k)}$ is negative. {Likewise
when $j$ is odd, $iF_{2k-1}$ has opposite signs at the  endpoints of
$I_{j,2k-1}$.}

 It follows that
$F_{2k-1}$ has an odd number of zeros within $I_{i,2k-1}$ as long as
it is contained in $[\pi/10, 9\pi/10]$.  Moreover, since the error
term $iR_{2k-1}(z_{\theta})$ is strictly less than the value of
$|2\sin(\theta(2k+1)/2)|$ outside of the intervals $I_{j,2k-1}$, we
can conclude that $F_{2k-1}(z_{\theta})$ has no zero in
$[\pi/10, 9\pi/10]\setminus \cup_j I_{j,2k-1}$. {The analogous
result holds for the zeros of $F_{2k+1}$ with respect to the
intervals $I_{j,2k+1}$.} As in the classical case, we can conclude
that within $[\pi/10, 9\pi/10]$, the zero intervals $I_{j,2k-1}$ for
$F_{2k-1}(z_{\theta})$ are separated by the zero intervals
$I_{j,2k+1}$ for $F_{2k+1}(z_{\theta})$.

As mentioned in Section 1, our method can be generalized to some other families of Eisenstein series,
in particular if the zeros of the family under consideration are expected to mainly locate on some
fixed hyperbolic line.  However, careful analysis is needed to handle the error terms in a more general setting.

\section{Zeros of Eisenstein series and special values of $L$-series}\label{L}

We now relate the Eisenstein series $E_{2k+1,\chi}(z)$
defined by \eqref{eq:E2k+1} to special values of a particular
$L$-series by defining
\begin{align}C:=\frac{4(-1)^k}{e_{2k}}=\frac{2(-1)^k}{(2k)!L(2k+1,\chi)}\left(\frac{\pi}{2}\right)^{2k+1}.\label{C}\end{align}
\noindent We have the following result in the style of \eqref{eq:2}.
\begin{theorem}\label{zeros} Let $k\geq 0$ be an integer.  Then
\begin{align*}\frac{2(-1)^k}{(2k)!L(2k+1,\chi)}\left(\frac{\pi}{2}\right)^{2k+1}=4(2k+1)-16\sum_{\tau\in\hp\backslash \Gamma(2)}\mathrm{ord}_\tau(E_{2k+1,\chi})\frac{1}{\lambda(\tau)}.\end{align*}
\end{theorem}

\begin{proof}
We start by building the function
\begin{align*}\widetilde{E}_{2k+1,\chi}(z):=\frac{E_{2k+1,\chi}(z)}{\Theta_2^{4k}(z)\Theta_3^2(z)}.\end{align*}
 By \cite[Lemma 2.2]{ly05}, $\Theta_3^2$ has the same character
as $E_{2k+1, \chi}$, thus $\widetilde{E}_{2k+1,\chi}$ is a modular function on $\G(2)$. By
examining the poles of $\widetilde{E}_{2k+1,\chi}$ we conclude that
it is a monic polynomial in $1/\lambda$ of degree $k$.  We write
\begin{align}\widetilde{E}_{2k+1,\chi}(z)=\left(\frac{1}{\lambda(z)}-\frac{1}{\lambda(\alpha_1)}\right)\left(\frac{1}{\lambda(z)}-\frac{1}{\lambda(\alpha_2)}\right)\cdots
\left(\frac{1}{\lambda(z)}-\frac{1}{\lambda(\alpha_k)}\right).\label{poly}\end{align}

Since $\Theta_2$ and $\Theta_3$ are holomorphic on $\hp$,
the $\alpha_j$ are exactly the zeros of $E_{2k+1,\chi}$ in
$\hp$, with multiplicity. Expanding the product in (\ref{poly}),
we have
\begin{align*}\widetilde{E}_{2k+1,\chi}(z)-\left(\frac{1}{\lambda(z)}\right)^k=-\left(\sum_{j=1}^k\frac{1}{\lambda(\alpha_j)}\right)\left(\frac{1}{\lambda(z)}\right)^{k-1}
+O\left(\frac{1}{\lambda(z)^{k-2}}\right).\end{align*} Comparing the
coefficients of the $q^{-(k-1)/2}$, we have
\begin{align*}16^{-k}(C-4-8k)=
\frac{-1}{16^{k-1}}\sum_{j=1}^k\frac{1}{\lambda(\alpha_j)}
=\frac{-1}{16^{k-1}}\sum_{\tau\in\hp\backslash
\Gamma(2)}\mathrm{ord}_\tau(E_{2k+1,\chi})\frac{1}{\lambda(\tau)}.\end{align*}
Solving for $C$ and using (\ref{C}) yields the desired result.
\end{proof}

\noindent We note that the set of zeros of $G_{2k+1}(z)=E_{2k+1,\chi}(z)\,|\,\gamma_0$ is given by
\begin{align*}\{\tau=\gamma_0^{-1}\tau_0:\tau_0\in\hp\backslash\Gamma(2),E_{2k+1}(\tau_0)=0\}.\end{align*}
Therefore, we have the following corollary.

\begin{cor}\label{cor:3.2} Let $k\geq 0$ be an integer.  Then
\begin{align*}\frac{2(-1)^k}{(2k)!L(2k+1,\chi)}\left(\frac{\pi}{2}\right)^{2k+1}=4(2k+1)-16\sum_{\tau\in\hp\backslash \Gamma(2)}\mathrm{ord}_\tau(G_{2k+1,\chi})\frac{1}{\lambda(\gamma_0\tau)}.\end{align*}
\end{cor}

\begin{proof} This follows directly from Theorem \ref{zeros} and the comment above.\end{proof}

We now consider Eisenstein series of even weight for $\Gamma(2)$.
With $S_k(x,y)$ defined as in \eqref{Sk} and $k\geq 1$ an integer,
we define
\begin{align*}
E_{2k}^{\pm}(z):=\frac{1}{2}(S_{2k}(0,1)\pm
S_{2k}(1,0)).\end{align*} These are modular forms of weight $2k$ for
$\Gamma(2)$ with trivial character $\chi_0$.  In \cite{ly05}, Long
and Yang show that $E^+_{2k}$ and $E^-_{2k}$ have the Fourier
expansions
\begin{align*}E_{2k}^{\pm}(z)=1\pm \frac{(2\pi i)^{2k}}{4^k\Gamma(2k)L(2k,\chi_0)}\sum_{r=1}^\infty\frac{r^{2k}q^{r/2}}{1+\mp (-1)^rq^{r/2}}.\end{align*}

\begin{theorem}\label{Eeven} Let $k\geq 1$ be an integer.  Then
\begin{align*}
\frac{(2\pi i)^{2k}}{4^k\Gamma(2k)L(2k,\chi_0)} =\pm 4(2k)+\mp
16\sum_{\tau\in\hp\backslash
\Gamma(2)}\mathrm{ord}_\tau(E_{2k}^{\pm})\frac{1}{\lambda(\tau)}.\end{align*}
\end{theorem}

\begin{proof}  To prove the first result, we set
\begin{align*}\widetilde{E}_{2k}^+(z):=\frac{E_{2k}^+(z)}{\Theta_2^{4k}}.\end{align*}
Then $\widetilde{E}_{2k}^+(z)$ is a modular function for $\Gamma(2)$
and is a monic polynomial in $1/\lambda$ of degree $k$.  Arguing as
in the proof of Theorem \ref{zeros}, we have
\begin{align*}\widetilde{E}_{2k}^+(z)-\left(\frac{1}{\lambda(z)}\right)^k=-\left(\sum_{\tau\in\hp\backslash \Gamma(2)}\mathrm{ord}_\tau(E_{2k}^+)\right)
\left(\frac{1}{\lambda(z)}\right)^{k-1}
+O\left(\frac{1}{\lambda(z)^{k-2}}\right).\end{align*} Equating the
coefficients of $q^{-(k-1)/2}$ on either side, we find that
\begin{align}16^{-k} \left (\frac{(2\pi i)^{2k}}{4^k\Gamma(2k)L(2k,\chi_0)}-8k \right )=-16^{-(k-1)}\sum_{\tau\in\hp\backslash \Gamma(2)}\mathrm{ord}_\tau(E_{2k}^+).\label{Eplus} \end{align}

Solving (\ref{Eplus}) yields the formula for $E_{2k}^+$. The one for
$E_{2k}^-$ is similarly obtained.
\end{proof}

\noindent Theorem \ref{Eeven} leads directly to the following surprising
identity.

\begin{cor} Let $k\geq 1$ be an integer.  Then
\begin{align*}k=\sum_{\tau\in\hp\backslash \Gamma(2)}(\mathrm{ord}_\tau(E_{2k}^+)+\mathrm{ord}_\tau(E_{2k}^-))\frac{1}{\lambda(\tau)}=\sum_{\tau\in\hp\backslash \Gamma(2)}\mathrm{ord}_\tau(E_{2k}^+E_{2k}^-)\frac{1}{\lambda(\tau)}.\end{align*}
\end{cor}

\providecommand{\bysame}{\leavevmode\hbox
to3em{\hrulefill}\thinspace}
\providecommand{\MR}{\relax\ifhmode\unskip\space\fi MR }
\providecommand{\MRhref}[2]{%
  \href{http://www.ams.org/mathscinet-getitem?mr=#1}{#2}
} \providecommand{\href}[2]{#2}

\end{document}